\documentclass[12pt]{article}
\usepackage{mathrsfs}
\usepackage{amsmath,amssymb}

\usepackage{amsfonts}
\usepackage{latexsym}
\usepackage{amsthm}
\usepackage{pictex}

\usepackage[T2A]{fontenc}     
\usepackage[cp1251]{inputenc} 
\usepackage{graphicx}

\def\lb{\label}

\newtheorem{theorem}{Theorem}[section]
\newtheorem{definition}{Definition}
\newtheorem{lemma}{Lemma}[section]

\begin{document}

\def\a{\alpha} \def\cA{{\cal A}} \def\bA{{\bf A}}  \def\mA{{\mathscr A}}
\def\b{\beta}  \def\cB{{\cal B}} \def\bB{{\bf B}}  \def\mB{{\mathscr B}}
\def\g{\gamma} \def\cC{{\cal C}} \def\bC{{\bf C}}  \def\mC{{\mathscr C}}
\def\G{\Gamma} \def\cD{{\cal D}} \def\bD{{\bf D}}  \def\mD{{\mathscr D}}
\def\d{\delta} \def\cE{{\cal E}} \def\bE{{\bf E}}  \def\mE{{\mathscr E}}
\def\D{\Delta} \def\cF{{\cal F}} \def\bF{{\bf F}}  \def\mF{{\mathscr F}}
\def\c{\chi}   \def\cG{{\cal G}} \def\bG{{\bf G}}  \def\mG{{\mathscr G}}
\def\z{\zeta}  \def\cH{{\cal H}} \def\bH{{\bf H}}  \def\mH{{\mathscr H}}
\def\e{\eta}   \def\cI{{\cal I}} \def\bI{{\bf I}}  \def\mI{{\mathscr I}}
\def\p{\psi}   \def\cJ{{\cal J}} \def\bJ{{\bf J}}  \def\mJ{{\mathscr J}}
\def\vT{\Theta}\def\cK{{\cal K}} \def\bK{{\bf K}}  \def\mK{{\mathscr K}}
\def\k{\kappa} \def\cL{{\cal L}} \def\bL{{\bf L}}  \def\mL{{\mathscr L}}
\def\l{\lambda}\def\cM{{\cal M}} \def\bM{{\bf M}}  \def\mM{{\mathscr M}}
\def\L{\Lambda}\def\cN{{\cal N}} \def\bN{{\bf N}}  \def\mN{{\mathscr N}}
\def\m{\mu}    \def\cO{{\cal O}} \def\bO{{\bf O}}  \def\mO{{\mathscr O}}
\def\n{\nu}    \def\cP{{\cal P}} \def\bP{{\bf P}}  \def\mP{{\mathscr P}}
\def\r{\rho}   \def\cQ{{\cal Q}} \def\bQ{{\bf Q}}  \def\mQ{{\mathscr Q}}
\def\s{\sigma} \def\cR{{\cal R}} \def\bR{{\bf R}}  \def\mR{{\mathscr R}}
\def\S{\Sigma} \def\cS{{\cal S}} \def\bS{{\bf S}}  \def\mS{{\mathscr S}}
\def\t{\tau}   \def\cT{{\cal T}} \def\bT{{\bf T}}  \def\mT{{\mathscr T}}
\def\f{\phi}   \def\cU{{\cal U}} \def\bU{{\bf U}}  \def\mU{{\mathscr U}}
\def\F{\Phi}   \def\cV{{\cal V}} \def\bV{{\bf V}}  \def\mV{{\mathscr V}}
\def\P{\Psi}   \def\cW{{\cal W}} \def\bW{{\bf W}}  \def\mW{{\mathscr W}}
\def\o{\omega} \def\cX{{\cal X}} \def\bX{{\bf X}}  \def\mX{{\mathscr X}}
\def\x{\xi}    \def\cY{{\cal Y}} \def\bY{{\bf Y}}  \def\mY{{\mathscr Y}}
\def\X{\Xi}    \def\cZ{{\cal Z}} \def\bZ{{\bf Z}}  \def\mZ{{\mathscr Z}}
\def\O{\Omega}
\def\ve{\varepsilon}
\def\vt{\vartheta}
\def\vp{\varphi}
\def\vk{\varkappa}

\def\ba{\breve{a}}
\def\bb{\breve{b}}

\newcommand{\gA}{\mathfrak{A}}
\newcommand{\gB}{\mathfrak{B}}
\newcommand{\gC}{\mathfrak{C}}
\newcommand{\gD}{\mathfrak{D}}
\newcommand{\gE}{\mathfrak{E}}
\newcommand{\gF}{\mathfrak{F}}
\newcommand{\gG}{\mathfrak{G}}
\newcommand{\gH}{\mathfrak{H}}
\newcommand{\gI}{\mathfrak{I}}
\newcommand{\gJ}{\mathfrak{J}}
\newcommand{\gK}{\mathfrak{K}}
\newcommand{\gL}{\mathfrak{L}}
\newcommand{\gM}{\mathfrak{M}}
\newcommand{\gN}{\mathfrak{N}}
\newcommand{\gO}{\mathfrak{O}}
\newcommand{\gP}{\mathfrak{P}}
\newcommand{\gR}{\mathfrak{R}}
\newcommand{\gS}{\mathfrak{S}}
\newcommand{\gT}{\mathfrak{T}}
\newcommand{\gU}{\mathfrak{U}}
\newcommand{\gV}{\mathfrak{V}}
\newcommand{\gW}{\mathfrak{W}}
\newcommand{\gX}{\mathfrak{X}}
\newcommand{\gY}{\mathfrak{Y}}
\newcommand{\gZ}{\mathfrak{Z}}

\def\mA{{\mathscr A}}
\def\mB{{\mathscr B}}
\def\mC{{\mathscr C}}
\def\mD{{\mathscr D}}
\def\mE{{\mathscr E}}
\def\mF{{\mathscr F}}
\def\mG{{\mathscr G}}
\def\mH{{\mathscr H}}
\def\mI{{\mathscr I}}
\def\mJ{{\mathscr J}}
\def\mK{{\mathscr K}}
\def\mL{{\mathscr L}}
\def\mM{{\mathscr M}}
\def\mN{{\mathscr N}}
\def\mO{{\mathscr O}}
\def\mP{{\mathscr P}}
\def\mQ{{\mathscr Q}}
\def\mR{{\mathscr R}}
\def\mS{{\mathscr S}}
\def\mT{{\mathscr T}}
\def\mU{{\mathscr U}}
\def\mV{{\mathscr V}}
\def\mW{{\mathscr W}}
\def\mX{{\mathscr X}}
\def\mY{{\mathscr Y}}
\def\mZ{{\mathscr Z}}

\def\Z{{\Bbb Z}}
\def\R{{\Bbb R}}
\def\C{{\Bbb C}}
\def\T{{\Bbb T}}
\def\N{{\Bbb N}}
\def\S{{\Bbb S}}
\def\H{{\Bbb H}}
\def\J{{\Bbb J}}
\def\dD{{\Bbb D}}
\def\W{{\Bbb W}}

\def\qqq{\qquad}
\def\qq{\quad}
\newcommand{\ma}{\begin{pmatrix}}
\newcommand{\am}{\end{pmatrix}}
\newcommand{\ca}{\begin{cases}}
\newcommand{\ac}{\end{cases}}
\let\ge\geqslant
\let\le\leqslant
\let\geq\geqslant
\let\leq\leqslant
\def\ma{\left(\begin{array}{cc}}
\def\am{\end{array}\right)}
\def\iint{\int\!\!\!\int}
\def\lt{\biggl}
\def\rt{\biggr}
\let\geq\geqslant
\let\leq\leqslant
\def\[{\begin{equation}}
\def\]{\end{equation}}
\def\wh{\widehat}
\def\wt{\widetilde}
\def\pa{\partial}
\def\sm{\setminus}
\def\es{\emptyset}
\def\no{\noindent}
\def\ol{\overline}
\def\iy{\infty}
\def\ev{\equiv}
\def\/{\over}
\def\ts{\times}
\def\os{\oplus}
\def\ss{\subset}
\def\h{\hat}
\def\Re{\mathop{\rm Re}\nolimits}
\def\Im{\mathop{\rm Im}\nolimits}
\def\supp{\mathop{\rm supp}\nolimits}
\def\sign{\mathop{\rm sign}\nolimits}
\def\Ran{\mathop{\rm Ran}\nolimits}
\def\Ker{\mathop{\rm Ker}\nolimits}
\def\Tr{\mathop{\rm Tr}\nolimits}
\def\const{\mathop{\rm const}\nolimits}
\def\dist{\mathop{\rm dist}\nolimits}
\def\diag{\mathop{\rm diag}\nolimits}
\def\Wr{\mathop{\rm Wr}\nolimits}
\def\BBox{\hspace{1mm}\vrule height6pt width5.5pt depth0pt \hspace{6pt}}

\def\Diag{\mathop{\rm Diag}\nolimits}


\def\Twelve{
\font\Tenmsa=msam10 scaled 1200 \font\Sevenmsa=msam7 scaled 1200
\font\Fivemsa=msam5 scaled 1200 \textfont\msbfam=\Tenmsb
\scriptfont\msbfam=\Sevenmsb \scriptscriptfont\msbfam=\Fivemsb

\font\Teneufm=eufm10 scaled 1200 \font\Seveneufm=eufm7 scaled 1200
\font\Fiveeufm=eufm5 scaled 1200
\textfont\eufmfam=\Teneufm \scriptfont\eufmfam=\Seveneufm
\scriptscriptfont\eufmfam=\Fiveeufm}

\def\Ten{
\textfont\msafam=\tenmsa \scriptfont\msafam=\sevenmsa
\scriptscriptfont\msafam=\fivemsa

\textfont\msbfam=\tenmsb \scriptfont\msbfam=\sevenmsb
\scriptscriptfont\msbfam=\fivemsb

\textfont\eufmfam=\teneufm \scriptfont\eufmfam=\seveneufm
\scriptscriptfont\eufmfam=\fiveeufm}

\title {Periodic Jacobi operator with  finitely supported perturbations}

\author{ Alexei Iantchenko
\begin{footnote}
{ Institute of Mathematics and Physics,
  Aberystwyth Univ.,
Penglais, Ceredigion,  SY23 3BZ, UK, email: aii@aber.ac.uk }
\end{footnote} \and
Evgeny Korotyaev
\begin{footnote}
{Sankt Petersburg, e-mail: korotyaev@gmail.com}
\end{footnote}
}

\maketitle

\begin{abstract}
 \no We describe the spectral properties of the  Jacobi operator
$(Hy)_n= a_{n-1}
y_{n-1}+a_{n}y_{n+1}+b_ny_n,$ $n\in\Z,$  with $
  a_n=a_n^0+ u_n,$
 $b_n= b_n^0+ v_n,$ where sequences  $a_n^0>0,$ $b_n^0\in\R$ are periodic with period $q$,
 and sequences $ u_n,$ $ v_n$ have
 compact support. In the case $ u_n\equiv 0$ we   obtain the asymptotics of the spectrum in the limit of small perturbations $ v_n.$
\end{abstract}


\section{Introduction}
\setcounter{equation}{0}
Let $H^0$ denote the $q$--periodic Jacobi matrix associated with the equation
\begin{equation}\lb{1e}
(H^0y)_n=a_{n-1}^0 y_{n-1}+a_{n}^0y_{n+1}+b_n^0y_n=\l
 y_n,\,\, (\l,n)\in\C\ts\Z
 \end{equation} where the sequences $a_n^0,b_n^0\in\R$ verify $a_{n}^0=a_{n+q}^0 >0,\,\,
 b_{n}^0=b_{n+q}^0\in\R,\,\,n\in\Z.$ We use also notation $a^0=(a_n^0)_{n\in\Z},$ $b^0=(b_n^0)_{n\in\Z}.$
In this paper we consider the spectral properties of a finitely supported perturbation  of $H^0$ defined as follows. Let $ u=( u_n)_{n\in\Z},$ $ v=( v_n)_{n\in\Z}$ be finite sequences satisfying
\begin{equation}\lb{potential}
 u_n, v_n\in\R,\,\, u_n=0,\,\, v_n=0\,\, \mbox{for}\,\, n<0\,\, \mbox{and}\,\, n>p,\,\, v_0, v_p\neq 0.
 \end{equation}
Let $H$ denote the infinite Jacobi matrix  associated with the equation
\begin{equation}\lb{pert}
(Hy)_n= a_{n-1}
y_{n-1}+ a_{n}y_{n+1}+ b_ny_n=\l y_n,\qq
  a_n=a_n^0+ u_n\qq
  b_n= b_n^0+ v_n,\,\, (\l,n)\in\C\ts\Z.
\end{equation}

For $a_n=1,$ $n\in \Z,$ the operator $H$ is the finite difference Schr{\"o}dinger operator with finitely supported potential.

 A lot of papers is devoted to the direct and inverse resonance problems for the Schr{\"o}dinger operator $-\frac{d^2}{dx^2}+q(x)$ on the line $\R$  with compactly supported perturbation (see \cite{S}, \cite{Fr}, \cite{Z}  and references given there).

The problem of resonances for the Schr{\"o}dinger with periodic plus compactly supported potential
$-\frac{d^2}{dx^2}+p(x)+q(x)$  is much less studied: \cite{F1}, \cite{KM}, \cite{K3}.

The inverse resonance problem is not yet solved. Finite-difference Schr{\"o}dinger and  Jacobi operators
express many similar features.
Spectral and scattering properties of infinite Jacobi matrices are much studied (see \cite{Mo}, \cite{DS1}, \cite{DS2} and references given there).
The inverse problem was solved for periodic Jacobi operators: \cite{P}.

 The inverse scattering problem  for asymptotically periodic coefficients was solved by Ag.~Kh. Khanmamedov: \cite{Kh1} (on the line, the russian versions are dated much earlier)  and I. Egorova, J. Michor, G. Teschl \cite{EMT} (on the line in case of quasi-periodic background).

The resonance problems are less studied (see M.Marletta and R.Weikard \cite{MW}). The inverse resonances problem was recently solved in the case of constant background \cite{K2}.
In the present article we consider the direct resonance problem in the case of periodic background. The results are applied in
in  \cite{IK1} to the inverse resonance problem.  In \cite{IK2} and \cite{IK3} we solved the direct and inverse resonance problems on the half-lattice and applied the results to the zigzag half-nanotube in magnetic field.

Let $\vp=(\vp_n(z))_{n\in\Z}$ and $ \vt=(\vt_n(z))_{n\in \Z} $ be
fundamental solutions for equation (\ref{1e}), under the conditions
$\vt_0=\vp_1=1$ and $\vt_1=\vp_0=0$. Let $\psi^\pm=\vt+m_\pm\vp$ be Floquet-Bloch functions (see Section \ref{ss-perJ}). Here $m_\pm$ are the Titchmarch-Weyl functions.

 Denote
$\Delta(\l)=2^{-1}(\vp_{q+1}+\vt_{q})$ the Lyapunov function.
Then it is known that  the zeros $\{E_j\}_{j=1}^{2q}$ of the polynomial $\Delta^2-1$ of
degree $2q$ can be enumerated as follows
$\l_0^+<\l_1^-\leq\l_1^+<\ldots<\l_{q-1}^-\leq\l_{q-1}^+<\l_q^-,$
where $\l_0^+=\l_0^-$ and $\l_q^-=\l_q^+.$  Then the spectrum of $H^0$ in $\ell^2(\Z)$  is absolutely continuous and consists
of $q$ zones $\sigma_n=[\l_{n-1}^+,\l_n^-],\,\,n=1,\ldots,q,$
separated by gaps $\gamma_n=(\l_n^-,\l_n^+).$
 In each gap there is one simple zero of
polynomials $\vp_q(\l),$ $\Delta (\l),$ $\vt_{q+1}(\l).$ Note that
$\Delta (\l_n^\pm)=(-1)^{q-n}.$

 We use the standard definition of the root: $\sqrt{1}=1$ and fix the branch of the function $\O(\l)=\sqrt{1-\Delta^2(\l)}$ on $\C$ by demanding
 $$
 \O(\l)=\sqrt{1-\D^2(\l)} <0,\qqq \mbox{for} \qq \l\in (\l_{q-1}^+, \l_{q}^-).
 $$
 Now we  introduce the two-sheeted Riemann surface $\Lambda$ of $\sqrt{\Delta^2(\l)-1}$
 obtained by joining the upper and lower rims of
two copies of the cut plane $\Gamma=\C\setminus \sigma_{\rm ac} (H^0)$ in
the usual (crosswise) way. The $k-$th gap on the first physical
sheet $\Lambda_+$ we will denote by $\gamma_k^+$ and the same gap
but on the second nonphysical sheet $\Lambda_-$  we will denote by
$\gamma_k^-$ and let $\gamma^{\rm c}_k$  be the union of
$\overline{\gamma^+_k}$ and $\overline{\gamma^-_k}$:
\begin{equation}\lb{union}
\gamma_k^{\rm c}=\overline{\gamma^+_k}\cup\overline{\gamma^-_k}.
\end{equation}

The finitely supported perturbation $( u, v)$  does not change the absolutely continuous spectrum:
 $\sigma_{\rm ac}(H)=\sigma_{\rm ac}(H^0)=\bigcup_{n=1}^q[\l_{n-1}^+,\l_n^-].$

Introduce the function $$\alpha(\l) = C\det ((H-\l)(H^0-\l)^{-1})=C\det
(I+(H-H^0)(H^0-\l)^{-1}),\qq C=\prod_{j=0}^p\frac{a_j}{ a_j},$$
 which is meromorphic on $\L$, see [F1].
Recall that $T=1/\alpha$ is the transmission coefficient in the
$S-$matrix for the pair $H,H^0$ (see Section \ref{s-pert}).   If $\alpha$ has some
poles, then they coincide with some $\l_k^\pm .$ It is well known
that if $\alpha(\l) =0$ for some zero $\l\in\Lambda_+,$ then $\l$ is
an eigenvalue of $H$  and $\l\in\cup\gamma_k^+.$ Note that there are
no eigenvalues on the spectrum $\sigma_{\rm
ac}(H^0)\subset\Lambda_+$ since $|\alpha (\l)|\geq 1$ on
$\sigma_{\rm ac} (H^0).$

With $\Omega(\l)=(1-\Delta^2(\l))^{1/2}$  we define the functions $A$, $J$  by
$$J(\l)=2\Omega(\l+i0)\Im\alpha(\l+i0),\qq A(\l)=\Re\alpha (\l+i0)
-1,\qq\mbox{for}\,\,\l\in\sigma\, (H^0)\subset \Lambda_+.
$$
These functions were introduced for the  Schr{\"o}dinger operator on $\R$ with periodic plus compactly supported potentials  by the second author in \cite{K1}.
   We  show that $A, J$ are polynomials on $\C$ and they are real on the real
line. Instead of the function $\alpha$ we consider the modified
function $\xi =2i\Omega\alpha$ on $\Lambda.$ We show that $\xi$
satisfies
\begin{equation}\lb{1.3}\xi =2i\Omega\alpha =2i\Omega (1+A)
-J\qq\mbox{on}\,\,\Lambda .
\end{equation}

Recall that $\Omega$
 is analytic on $\Lambda$ and
$\Omega = 0$ for some $\l\in\Lambda$ iff $\l =\l_k^-$ or $\l
=\l_k^+$ for some $k\geq
 0.$ Then the function $\xi$ is analytic on $\Lambda$ and has branch points $\l_n^\pm,$ if $\gamma_n\neq \emptyset.$  The zeros of $\xi$  define the eigenvalues and
the resonances.  Define the set
$$\Lambda_0=\{\l\in\Lambda:\,\,\l
 =\l_k^+\in\Lambda_+\,\,\mbox{and}\,\, \l
 =\l_k^+\in\Lambda_-,\,\, \gamma_k=\emptyset\}\subset\Lambda.$$
  In fact with each $\gamma_k =\emptyset$ we associate two points
  $\l_k^+\in\Lambda_+$ and $\l_k^+\in\Lambda_-$ from the
set $\Lambda_0.$  If each gap of $H^0$ is not empty, then
$\Lambda_0=\emptyset.$

\begin{definition}\lb{states} Each zero of $\xi$ on $\Lambda\setminus\Lambda_0$
is a state of $H.$\\
1) A state $\l\in\Lambda_+$ is a bound state.\\
2) A state $\l\in\Lambda_-$ is a resonance.\\
3) A state $\l=\l_k^\pm,$ $k=1,\ldots,q$ is a virtual state.\\
A resonance $\l\in \cup\gamma_k^-\subset\Lambda_-$ is an anti-bound
state.
\end{definition}
This definition is motivated by the representation (\ref{resolvent}) of the resolvent.
It is known that the gaps $\gamma_k=\emptyset$ do not give
contribution to the states. Recall that $S-$matrix for $H,H^0$ is
meromorphic on $\Lambda,$ but it is analytic at the points from
$\Lambda_0$ (see [F1]). Roughly speaking there is no difference
between the points from $\Lambda_0$  and other points inside the
spectrum of $H^0.$


We prove the following result.
\begin{theorem}[Total number of states]\label{Th-number-of-states} Let $a^0,b^0$ be $q$--periodic sequences as in (\ref{1e}) and  $ u, v$  be a finitely supported perturbations satisfying (\ref{potential}). If $ u_p\neq0$ then
 $H$ has in total $2p+2q-1$
states. If $ u_p=0$ then the number of states is    $2p+2q-2.$
\end{theorem}
Theorem \ref{Th-number-of-states} follows from the properties of function $\xi(\l)$
formulated in Theorem \ref{Th-asymptotics-number-of-states}.

\begin{theorem} Let $a^0,b^0$ be $q$--periodic sequences as in (\ref{1e}) and  $ u, v$  be  finitely supported perturbations satisfying (\ref{potential}). Then $\xi$ satisfies (\ref{1.3}) and the following properties hold.\\
i) There exists an even number of states (counted with multiplicities) on each set $\gamma_k^{\rm c}\neq\emptyset,$ $k=1,\ldots,q-1,$ where $\gamma_k^{\rm c}$ is a union of the physical gap $\overline{\gamma}^+_k$ and non-physical gap $\overline{\gamma}^-_k$ (see (\ref{union})).\\
ii) Let $\l\in\gamma_k^+$ be a bound state for some $k=0,\ldots,q,$
i.e. $\xi(\l)=0.$ Let $\tilde{\l}\in\gamma_k^-\subset\Lambda_-$ be the same number but on the second sheet $\Lambda_+.$ Then $\tilde{\l}\in\gamma_k^-$ is not an antibound state i.e. $\xi(\tilde{\l})\neq 0.$\\
iii) Let $N$ be the total number of bound states. Then $N+q$ is even.
\end{theorem}
The properies i),ii) follows from Lemma \ref{Lemma_3.3_Korotyaev}.
The property iii) is proved in \cite{IK1}.

For each $k=1,\ldots,q-1$ there exists a unique point $\alpha_k\in  [\l_k^-,\l_k^+]$ such that
\begin{equation}\lb{1.4}
(-1)^{q-k}\Delta(\alpha_k)=\max_{\l\in [\l_k^-,\l_k^+]}|\Delta(\l)|=\cosh h_k
\end{equation}
for some $h_k\geq  0.$

Let $\mu_k$ , $k > 1$ be the Dirichlet spectrum of the equation
$H^0y=\l y_n$ on the interval $[0, q]$ with the boundary condition
$y_0 = y_q = 0.$ They are zeros of $\vp_q(\l).$ It is well known
that each $\mu_k\in [\l_k^- ,\l_k^+ ],$ $k=1,\ldots, q-1
> 1.$
\begin{theorem}[Small perturbations]\lb{Th_representation} Let $t> 0.$ Suppose $ u_j=0,$ and write $t v_j$ instead of $ v_j.$ Let $\gamma_n\equiv (\l_n^-,\l_n^+)\neq \emptyset.$ Then there exists $t_0>0$ small enough  such that for all $t\in (0,t_0)$ the following statements hold true.\\
 i)  In each $\overline{\gamma}_k\neq\emptyset,$ $k=1,\ldots,q-1,$ there are exactly two simple real states  $\l^\pm_k(t)\in\gamma_k^{\rm c}$ such that $\l_k^-\leq \l_k^-(t) < \alpha_k <\l_k^+(t)\leq\l_k^+.$
  Moreover, if $\l_0$ is one of $\l_n^\pm(t)$ and satisfies $(-1)^{q-n+1}J(\l_0) >0$ (or $(-1)^{q-n+1}J(\l_0) <0$ or $J(\l_0)=0$), then $\l_0$ is a bound state (or an anti-bound state or a virtual state).\\
  ii)   The following asymptotics hold true
$$\l^\pm_n(t)=\l_n^\pm+t^2\frac{J_1^2(\l_n^\pm)}{4(\Delta^2)'(\l_n^\pm)}+{\mathcal O}(t^3),$$
where $J(\l_n^\pm)=tJ_1(\l_n^\pm)+{\mathcal O}(t^2)$ and, if $\l_n^\pm\neq\mu_n$ then
\begin{align*}
J_1(\l_n^\pm)=-\frac{\vp_q(\l_n^\pm)}{a_0^0}\sum_{k=0}^p v_k|\psi_k^+(\l_n^\pm)|^2,
\end{align*} otherwise if $\l_n^+=\mu_n$ or $\l_n^-=\mu_n$ then
$$J_1(\mu_n)=\frac{\vt_{q+1}(\mu_n)}{a_0^0}\sum_{k=1}^{p} v_k(\vp_k(\mu_n))^2.$$
\end{theorem}
{\bf Remark.} Suppose $ v_k >0,$ $k=1,\ldots, p$ and $\l_n^\pm\neq\mu_n.$ Then for $n=1,\ldots, q-1,$ we have $(-1)^{q-n+1}J_1(\l_n^-) >0$ and  $\l_n^-(t)$ is bound state and $\l_n^+(t)$ is antibound state.  If $\l_n^-=\mu_n$ then  we get that $\l_n^-(t)$ is antibound state and $\l_n^+(t)$ is bound state. Similar in other cases. We have that if all $ u_k,$ $k=0,\ldots,p,$ have constant sign then for $t$ small enough there is exactly one bound state in each open gap $\gamma_k^+\neq\emptyset,$ $k=1,\ldots,q-1.$

If in the contrary the sign of $ u_k$ is not constant then there can be open gaps $\gamma_n$ with exactly two bound states or exactly two antibound states, namely if $$\sign \sum_{k=0}^p v_k|\psi_k^+(\l_n^-)|^2=\sign \sum_{k=0}^p v_k|\psi_k^+(\l_n^+)|^2,\qq\mbox{if}\,\, \l_n^\pm\neq\mu_n,$$
and similar in the other cases.

\section{Periodic Jacobi matrices.}\label{ss-perJ}

We recall some known properties of the $q-$periodic Jacobi matrix
\[
\lb{21}
 H^0=\left(\begin{array}{cccccccc}
  ...& ...  & ...  & ...     &...         & ...        &... \\
 ...&a_0^0  & b_1^0  & a_1^0     &0         & 0        &... \\
 ...&0  & a_1^0& b_2^0       &a_2^0& 0        &... \\
 ...&0  & 0  &a_2^0& b_3^0         & b_3^0     &... \\
 ...&0  & 0  & 0        &a_3^0      & b_4^0       &... \\
 ...&0  & 0  & 0        &0         &a_4^0&... \\
 ...&...  & ...& ...      &...       &...       &... \\
 \end{array}\right),\qq a_{n}^0=a_{n+q}^0 >0,\,\,
 b_{n}^0=b_{n+q}^0\in\R,\,\,n\in\Z,
\]
 associated with the equation (\ref{1e}): $
 (H^0y)_n=a_{n-1}^0 y_{n-1}+a_{n}^0y_{n+1}+b_n^0y_n=\l y_n,$
$(\l,n)\in\C\ts\Z .$

As before we denote $\Delta(\l)=2^{-1}(\vp_{q+1}+\vt_{q})$ the Lyapunov function, where
 $\vp=(\vp_n(\l))_{n\in\Z}$ and $ \vt=(\vt_n(\l))_{n\in \Z} $ are the
fundamental solutions for equation (\ref{1e}) satisfying
$\vt_0=\vp_1=1$ and $\vt_1=\vp_0=0$.

Let $$ A=\prod_{j=1}^qa_j^0>0,\,\,B=\sum_{j=1}^qb_j^0.$$

We have
\begin{equation}\lb{2.1}
\vt_n (\l)=-\frac{a_0^0}{a_1^0\ldots a_{n-1}^0}\l^{n-2} +{\mathcal
O}(\l^{n-3}),\qq\vp_n(\l)=\frac{\l^{n-1}}{a_1^0\ldots a_{n-1}^0}
+{\mathcal O}(\l^{n-2}).
\end{equation}
We denote the zeros of
 ${\vp_q}$ resp. ${\vt_{q+1}}$  by $\mu_n\in\gamma_n,$ resp. $\nu_n\in\gamma_n,$ $n=1,q-1$ (Dirichlet or Neumann eigenvalues).
Then
$$\vp_q=\frac{a_0^0}{A}\prod_{j=1}^{q-1}(\l-\mu_j),\qq \vt_{q+1}=\frac{-a_0^0}{A}\prod_{j=1}^{q-1}(\l-\nu_j).$$
 The zeros $\{E_j\}_{j=0}^{2q-1}$ of the polynomial $\Delta^2-1$
of degree $2q$ can be enumerated as follows
\begin{equation*}
E_0\equiv\l_0^+<E_1\equiv\l_1^-\leq
E_2\equiv\l_1^+<\ldots<E_{2q-3}\equiv \l_{q-1}^-\leq
E_{2q-2}\equiv\l_{q-1}^+<E_{2q-1}\equiv\l_q^-,\end{equation*} where
$\l_0^+=\l_0^-$ and $\l_q^-=\l_q^+.$  Then the
  spectrum of
$H^0$ on $\ell^2(-\infty,+\infty)$ is absolutely continuous and consists
of $q$ zones $\sigma_n=[\l_{n-1}^+,\l_n^-],\,\,n=1,\ldots,q$
separated by gaps $\gamma_n=(\l_n^-,\l_n^+).$ We denote $\gamma_0=(-\infty, \l_0^+),$ $\gamma_q=( \l_q^+,+\infty),$ the infinite gaps.
 In each gap $\gamma_n,$ $n=1,\ldots,q-1,$ there is one simple zero of
polynomials $\vp_q(\l),$ $\Delta (\l),$ $\vt_{q+1}(\l).$ Note that
$\Delta (\l_n^\pm)=(-1)^{q-n}.$

Let $\Gamma$ be complex $\l$-plane with cuts along segments
$\sigma_n,$ $n=1,2,\ldots,q.$  $\Gamma$ can be identified with $\Lambda^+$ and on $\Gamma$ we omit the index $+.$ On the plane $\Gamma$ consider the
function $$z=z(\l)=\Delta (\l) +\sqrt{\Delta^2(\l)-1},$$ fixing the
branch by the condition $\sqrt{\Delta^2(\l)-1} <0$ for $\l >
\l_q^-$ (in accordance with Introduction). Then
\begin{equation}\lb{Dfactor}
\sqrt{\Delta^2(\l)-1}=\frac{-1}{2A}\prod_{j=0}^{2q-1}\sqrt{\l-E_j}.
\end{equation}
 Then function $z=z(\l)$ is continuous up
to the boundary $\partial\Gamma$  and has the properties:\\
$|z| <1\,\,\mbox{for}\,\,\l\in\Gamma,$ and
$|z|=1\,\,\mbox{for}\,\,\l\in\partial\Gamma.$
Moreover (Teschl page 116 (7.12)) for $\l\in\Lambda_+,$
$$z^{\pm 1}\equiv\xi^\pm(\l)=(2\Delta(\l))^\mp\left(1+{\mathcal
O}\left(\l^{-2q}\right)\right)=\left(\frac{A}{\l^q}\right)^{\pm
1}\left( 1\pm \frac{B}{\l}+{\mathcal
O}\left(\frac{1}{\l^2}\right)\right).$$
We have
$$2i\Omega=z-z^{-1}=\mp\frac{\l^q}{A}\left( 1\mp
\frac{B}{\l}+{\mathcal O}\left(\frac{1}{\l^2}\right)\right),\qq \l\in\Lambda_\pm.$$

Then $z(\l)=\xi_+(\l)=e^{iq\varkappa(\l)}$ is the  first Floquet
multiplier and $\varkappa(\l)$ is quasimomentum. The second Floquet
multiplier is then $\xi_-(\l)=\overline{z(\l)}.$ We
denote also $e^{i\varkappa(\l)}=\omega,$ so
$z(\l)=\xi_+(\l)=\omega^q.$

For each $k=1,\ldots,q-1$ there exists an unique point $\alpha_k\in  [\l_k^-,\l_k^+]$ such that $
(-1)^{k-q}\Delta(\alpha_k)=\max_{\l\in [\l_k^-,\l_k^+]}|\Delta(\l)|=\cosh h_k.
$

On each $\gamma_k^+,$ $k=0,1,\ldots,q,$ the quasi-momentum $\varkappa(\l)$
 has  constant real part $\Re(\vk)= \frac{q-k}{q}\pi,$ $\vk(\l_k^-)=\vk(\l_k^+)=\frac{q-k}{q}\pi,$ and positive $\Im(\vk).$ Moreover, as $\l$ increases from $\l_k^-$ to $\l_k^+,$ the  imaginary part  $\Im(\vk)\equiv v(\l)$ starts by increasing from $0$ to $h_k=\frac{1}{i}(\vk(\alpha_k) -\frac{q-k}{q}\pi)$, then decreases from $h_k$ to $0.$ Then
 \[\lb{isin}\sqrt{\Delta^2(\l)-1}=
i\sin q\varkappa=2^{-1}e^{i(q-k)\pi}(e^{-qv}-e^{qv})=   -(-1)^{q-k}\sinh qv,\]
and also $\sinh qv=-2^{-1}(\omega^q-\omega^{-q})>0.$

Equation (\ref{1e}) has two solutions (Bloch functions) $\psi_n^\pm=\psi_n^\pm(\l)$ which satisfy
$\psi_{kq}^\pm=\xi_\pm^k,$ $k\in\Z,$ and at the end points of the gaps we have $|\psi^\pm_{kq}(\l_n^\pm)|=1.$
As for any $\l\in\Gamma$ we have $\psi^\pm\in \ell^2(0,\pm\infty),$
then functions $\psi^\pm(\l)$ are the Floquet solutions for
(\ref{1e}): \begin{align}&\psi_n^\pm
(\l)=\vt_n(\l)+m_\pm(\l)\vp_n(\l),\lb{Jost}\\
&m_\pm(\l)=
\frac{\phi\pm\sqrt{\Delta^2(\l)-1}}{\vp_q},\,\,\phi=\frac{\vp_{q+1}-\vt_q}{2},\,\,\l\in\Gamma.\lb{TitchWeyl}
\end{align}
Thus \begin{equation}\lb{mm} m_+m_-=-\frac{\vt_{q+1}}{\vp_q}.
\end{equation}
 This equality
considered at zeros of polynomial  $\vp_q(\l)$ shows
that  if  $\mu_n,$  $n=1,\ldots,q-1,$  is not a
virtual state then we have either (i) $m_+$ has simple pole at
$\mu_n,$ $m_-$ is regular or (ii) $m_-$ has simple pole at $\mu_n,$
$m_+$ is regular. Thus
one of the solutions $\psi_n^\pm(\l)$ is regular, the other has
simple poles, one in each
 finite gap $\gamma_n,$ $n=1,\ldots,q-1$.

We have also
\begin{equation}\label{Fact7}
m_+=\frac{z-\vt_q}{\vp_q},\,\,m_-=\overline{m_+},\qq\phi^2+1-\Delta^2=1-\varphi_{q+1}\vartheta_q=-\vartheta_{q+1}\varphi_q.
\end{equation}

\section{Scattering theory}\lb{s-pert}
For a finitely supported sequences $ u, v$ satisfying (\ref{potential}) we  consider the infinite Jacobi matrix $H$ associated with
 the equation (\ref{pert}):
$$
 (Hy)_n= a_{n-1} y_{n-1}+ a_{n}y_{n+1}+ b_ny_n=\l y_n,\qq
  a_n=a_n^0+ u_n\qq
  b_n= b_n^0+ v_n\qq
$$
$(\l,n)\in\C\ts\Z .$

Let $f^\pm$ be the Jost functions
$  f_n^-=\p_n^-,\,\, \mbox{for}\,\, n \leq 0;\,\,\,\,f_n^+=\p_n^+,\,\,\mbox{for}\,\,
n > p. $ Let
$\{\phi_n,\psi_n\}= a_n(\phi_n\psi_{n+1}-\phi_{n+1}\psi_n)$
denote the  Wronskian.

We have
$f_n^+=\tilde{\vt}_n +m_+\tilde{\vp}_n,$ $f_n^-=\hat{\vt}_n
+m_-\hat{\vp}_n,$ $m_-(\l)=\overline{m_+(\overline{\l})},$ where
$$\hat{\vt}_n,\hat{\vp}^-_n= \vt_n,\vp_n\,\,\mbox{for}\,\,n\leq
0\qq\mbox{and}\qq \tilde{\vt}_n,\tilde{\vp}_n=
\vt_n,\vp_n\,\,\mbox{for}\,\,n>p.$$

The Jost
functions are linearly independent and we have
\begin{equation}\lb{f-relations}
f^\pm_n(\l)=\alpha
(\l)\overline{f^\mp_n(\l)}+\beta_\mp(\l)f^\mp_n(\l),\,\,\mbox{for}\,\,\l\in\sigma_{\rm
ac}(H^0),
\end{equation}
where \begin{align*} \alpha (\l)&=\frac{\{f^\mp(\l),
f^\pm(\l)\}}{\{f^\mp(\l),
\overline{f^\mp(\l)}\}}=\frac{\varphi_q(\l)}{a_0^02i\sin q\varkappa
(\l)}\{f^-(\l), f^+(\l)\},\\
\beta_\pm (\l)&=\frac{\{f^\mp(\l),
\overline{f^\pm(\l)}\}}{\{f^\pm(\l),
\overline{f^\pm(\l)}\}}=\mp\frac{\varphi_q(\l)}{a_0^02i\sin q\varkappa
(\l)}\{f^\mp(\l),\overline{ f^\pm(\l)}\}
\end{align*}
 We denote $$s=
\{f^+(\l),\overline{ f^-(\l)}\},\qq w=\{f^-(\l), f^+(\l)\},$$
$$\Omega (\l)=(1-\Delta^2(\l))^{1/2}=\sin q\vk (\l),\qq\xi=2i\Omega\alpha,\qq \cF=\xi\xi^*=\frac{\varphi_q^2(\l)}{(a_0^0)^2}ww^*,\qq\cS=\frac{\varphi_q^2(\l)}{(a_0^0)^2}ss^*.$$

 \begin{equation}\lb{alpha}\alpha=\frac{\varphi_q(\l)}{a_0^02i\Omega} w(\l),\qq \beta_+=-\frac{\varphi_q(\l)}{a_0^02i\Omega} \overline{s}(\l),\qq \beta_-=\frac{\varphi_q(\l)}{a_0^02i\Omega} s(\l),\qq\xi= \frac{\varphi_q(\l)}{a_0^0} w(\l)\end{equation}

 Then $\overline{s}=
-\{f^-(\l),\overline{ f^+(\l)}\}$ and we have
\begin{align}&|\alpha(\lambda)|^2 =1
+|\beta_\pm(\l)|^2,\,\,\mbox{and}\,\,\overline{\beta_\pm(\l)}=-\beta_\mp(\l),\lb{3}\\
&ww^*=\frac{(a_0^0)^24\Omega^2}{\varphi_q^2(\l)}+ss^*,\qq \cF=4\Omega^2+\cS,\qq \Omega^2 (\l)=1-\Delta^2(\l)
\end{align}

 We define the scattering matrix
\begin{equation}\lb{Smatrix}
S(\l)=\left(
          \begin{array}{cc}
            T(\l) & R_-(\l) \\
            R_+(\l) & T(\l) \\
          \end{array}
        \right),\,\,\l\in\sigma (H^0),
\end{equation} for the pair $(H,H^0),$ where
$$T(\l)=\frac{1}{\alpha(\l)},\,\,R_\pm(\l)=\frac{\beta_\pm(\l)}{\alpha
(\l)}=\frac{\mp\{f^\mp(\l),\overline{ f^\pm(\l)}\}}{\{f^-(\l),
f^+(\l)\}}.$$
We have also $R_+=\overline{s}/w,$ and $R_-=s/w.$

 The matrix $S(\l)$ is unitary:
$|T(\l)|^2+|R_\pm(\l)|^2=1,$
$T(\l)\overline{R_+(\l)}=-\overline{T(\l)} R_-(\l),$ extends to $\Lambda$ as a meromorphic function. The quantities
$T$ and $R_\pm$ are the transmission and the reflection coefficients
respectively:
$$T(\l)f^\pm_n(\l)=\left\{\begin{array}{lr}
                            T(\l) \psi^\pm_n(\l),& n\rightarrow \pm\infty \\
                            \psi_n^\pm(\l) +R_\mp(\l)\psi_n^\mp(\l),
                            & n\rightarrow \mp\infty
                          \end{array}\right.,\,\,\,\,\l\in\sigma
                          (H^0).
$$

The determinant of the scattering matrix is given by
$$\det S(\l)=T^2-R_+R_-=\frac{1}{\alpha^2} +\frac{|\beta|^2}{\alpha^2}=\frac{|\alpha|^2}{\alpha^2}=\frac{\overline{\alpha (\l)}}{\alpha (\l)}.$$

We have also $$R_+=\frac{\overline{s}}{w},\qq R_-=\frac{s}{w}.$$

\begin{lemma}\lb{l-3.2} The following identities hold true.
\begin{align}
&\xi(\l)=2i\sin q\vk (1+A)-J,\qq A= \frac12\left[
\left(\frac{ a_0}{a_0^0}\tilde{\vp}_1-\vp_1\right)+\frac{ v_0}{a_0^0}\tilde{\vp}_0+(\tilde{\vt}_0-\vt_0)
\right],\nonumber\\
&J=-\left[\frac{ a_0}{a_0^0}\vp_q\tilde{\vt}_1+\phi(\frac{ a_0}{a_0^0}\tilde{\vp}_1-\tilde{\vt}_0)+
\frac{ v_0}{a_0^0}(\vp_q\tilde{\vt}_0+\phi
\tilde{\vp}_0)+\vt_{q+1}\tilde{\vp}_0 \right],\nonumber\\
&\xi(\l)=2(-1)^{q-n+1}\sinh v(\l)(1+A) -J(\l),\qq \l\in\gamma_n^\pm\neq\emptyset, \lb{3.10}
\end{align}
where $v=\Im\vk$ and $\pm v(\l) >0$ for $\l\in\gamma_n^\pm.$

\end{lemma}
{\bf Proof.}
Let us calculate $\alpha.$ We have $\psi_0^-=1,$ $\psi_1^-=m_-,$
then using
$$y_{n-1}=\frac{(\l-b_n^0)y_n-a_n^0y_{n+1}}{a_{n-1}^0}$$ we get
\begin{align*}
&\psi_{-1}^-=\frac{(\l -b_0^0)\cdot\psi_0^-
-a_0^0\psi_1^-}{a_{-1}^0}=\frac{(\l -b_0^0) -a_0^0m_-}{a_{-1}^0}.
\end{align*}
Now applying
$$y_{n+1}=\frac{(\l- b_n)y_n- a_{n-1}y_{n-1}}{ a_n}.$$
 Now,
$f_0^-=\psi_0^-=1,$ $f_{-1}^-=\psi_{-1}^-.$ Then
$$f_1^-=\frac{(\l - b_0)\cdot f_0^-
-a_{-1}^0f_{-1}^-}{ a_0}=\frac{(\l- b_0)-(\l
-b_0^0)\cdot\psi_0^- +a_0^0\psi_1^-}{ a_0}=\frac{a_0^0m_-
- v_0}{ a_0} .$$
Recall
$$m_+=\frac{\omega^q-\vt_q}{\vp_q},\,\,m_-=\overline{m_+}.$$

 For $|\omega|=1,$ $\omega^2\neq 1,$ \begin{align*}&
 w=\{f_n^-,f_n^+\}=\const=\{f_0^-,f_0^+\}= a_0(f_0^-f_1^+-f_1^-f_0^+)= a_0f_1^++( v_0-a_0^0m_-)
 f_0^+,\\
  &s=\{f_n^+,\tilde{f}_n^-\}=\const=\{f_0^+,\tilde{f}_0^-\}=
  a_0(f_0^+\tilde{f}_1^--f_1^+\tilde{f}_0^-)= (a_0^0\tilde{m}_-- v_0)
 f_0^+
 - a_0f_1^+.
  \end{align*}

Now using (\ref{alpha}) we get $$\xi(\l)=2i\sin q\vk\cdot\Re\alpha-2\sin
q\vk\cdot\Im\alpha,\,\,\l\in\sigma_{\rm ac}(H^0).$$ Then we have
$\cF=\xi\xi^*=4(1-\Delta^2)(\Re\alpha)^2+(2\sin
q\vk\cdot\Im\alpha)^2.$ Denote $A=\Re\alpha -1$ and $J=2\sin
q\vk\cdot\Im\alpha.$ Hence
\begin{equation}\label{xi}
\xi(\l)=2i\sin q\vk (1+A)-J.
\end{equation}

We have
\begin{align*}
&J=-\left[\frac{ a_0}{a_0^0}\vp_q\tilde{\vt}_1+\frac{ a_0}{a_0^0}\phi\tilde{\vp}_1+\frac{ v_0}{a_0^0}\vp_q\tilde{\vt}_0+\frac{ v_0}{a_0^0}\phi
\tilde{\vp}_0-\phi\tilde{\vt}_0+\vt_{q+1}\tilde{\vp}_0.
\right]\\
&=-\left[\frac{ a_0}{a_0^0}\vp_q\tilde{\vt}_1+\phi(\frac{ a_0}{a_0^0}\tilde{\vp}_1-\tilde{\vt}_0)+
\frac{ v_0}{a_0^0}(\vp_q\tilde{\vt}_0+\phi
\tilde{\vp}_0)+\vt_{q+1}\tilde{\vp}_0 \right].
\end{align*}
As for $\l\in\sigma_{ac}(H^0),$ $\Im m_+=\vp_q^{-1}\sin q\kappa,$ we
get
\begin{equation}\label{A}
\Re\alpha=\frac12\left[\frac{ a_0}{a_0^0}\tilde{\vp}_1+\frac{ v_0}{a_0^0}\tilde{\vp}_0+\tilde{\vt}_0
\right],\qq A= \frac12\left[
\left(\frac{ a_0}{a_0^0}\tilde{\vp}_1-\vp_1\right)+\frac{ v_0}{a_0^0}\tilde{\vp}_0+(\tilde{\vt}_0-\vt_0)
\right]. \end{equation}
Identity (\ref{3.10}) follows from (\ref{isin}).
\qed

Lemma \ref{l-3.2} shows that $\xi$ is analytic on $\Lambda.$ We define the states by Definition \ref{states}.

{\bf Resolvent.}   The kernel of the resolvent of $H$ is
\begin{equation}\lb{resolvent}
R(n,m)=\langle e_n,(J-\lambda)^{-1} e_m\rangle=-
\frac{f_n^- f_m^+}{\left\{f^-,f^+\right\}}=-
\frac{F_{n,m}}{a_0^0\xi},\,\,n <m,
\end{equation}
 where
$e_n=(\delta_{n,j})_{j\in\Z},$ $F_{n,m}=\vp_q f_n^- f_m^+.$

The function $R(n,m)$ is meromorphic on $\Lambda$ for each
$n,m\in\Z.$

We have $$F_{n,m}=\vp_q\tilde{\vartheta}_n
\tilde{\vartheta}_m+(\phi+i\sin q\vk)
\tilde{\vartheta}_n\tilde{\varphi}_m+(\phi-i\sin q\vk)
\tilde{\vp}_n\tilde{\vt}_m-\vartheta_{q+1}\tilde{\varphi}_n\tilde{\varphi}_m.$$
The zeros of $\xi$ define the bound states and resonances as
$F_{n,m}$ is locally bounded. This motivates Definition \ref{states}.

We get $\cF =4(1-\Delta^2)(1+A)^2+J^2,$ where $J,A$ are
polynomials, and $A,J\equiv 0$ if $ u_j, v_j\equiv 0$
for all $j\in\Z.$

Using (\ref{3}) and $4(1-\Delta^2)=(2\sin q\vk)^2$
we have also $\cF =4(1-\Delta^2)+\cS(\l),$ where
$\cS(\l)=\frac{\vp^2_q(\l)}{(a_0^0)^2}s(\l)s^*(\l).$ Valid for all $\l\in\C.$

\begin{lemma}\lb{Lemma_3.3_Korotyaev}
1) We have $\cF =4(1-\Delta^2)(1+A)^2+J^2=4(1-\Delta^2)+\cS(\l)$ is polynomial and $\cF(\l)>0 $ and $\cS(\l)\geq 0$ on each interval
$(\l_{k-1}^+,\l_k^-),$ $k=1,\ldots, q.$ The function $\cF$ has even
number of zeros on each interval $[\l_k^-,\l_k^+],$ $k=1,\ldots,q-1$
and $\cF$ has only simple zeros
at $\l_k^\pm,$ $\gamma_k\neq\emptyset.$\\
2) If $\gamma_k=(\l_k^-,\l_k^+)=\emptyset$ for some $k\neq 0,$ then
each $f^\pm_n(\cdot),$ $n\in\Z,$ is analytic in some disk
$B(\l_k^+,\epsilon),$ $\epsilon >0.$ Moreover, $\mu_k=\l_k^\pm$ is a
double zero of $\cF$ and $\l_k^\pm$ is not a state of $J.$\\
3) Let $\l\in\gamma_k^+$ be a bound state for some $k=0,\ldots,q,$
i.e. $\xi(\l)=0.$ Then $\l\in\gamma_k^-$ is not an antibound state
and $\xi(\l)\neq 0.$\\
4) $\l\in\C$ is a zero of $\cF$ iff $\l\in\Lambda$ is a zero of
$\xi$ with the same multiplicity.
\end{lemma}
{\bf Proof.} 1) As $\cF =4(1-\Delta^2)(1+A)^2+J^2,$ it is clear that
$\cF(\l)>0 $ and $\cS(\l)\geq 0$ on each interval
$(\l_{k-1}^+,\l_k^-),$ $k=1,\ldots, q.$ Due to $\cF(\l_k^\pm)\geq
0,$ we get that $\cF$ has  even number of zeros in each interval
$[\l_k^-,\l_k^+],$ $k=1,\ldots,q-1.$  We have $\cF=\cF_0+\cS,$
$\cF_0=4(1-\Delta^2).$ We consider the case $\l=\l_k^+,$ the proof for
$\l=\l_k^-$ is similar. If $\cF(\l)=0,$ then we get $\cS(\l)=0,$
since $\Delta^2(\l_k^+)=1.$ Moreover,
$\cF_0'(\l_k^+)=-2\Delta(\l_k^+)\Delta'(\l_k^+) >0$ and
$\cS'(\l_k^+)\geq 0,$ which gives that $\l=\l_k^+$ is a simple zero
of $\cF.$

2) Suppose $\mu_k=\l_k^-=\l_k^+.$  The function $m_\pm$ is analytic
on $\Gamma=\C\setminus \cup\overline{\gamma}_k$ and hence each
$f^\pm_n(\cdot),$ $n\in\Z,$ is analytic in $\Gamma.$ Moreover
(\ref{Fact7}) yields $\phi(\mu_k)=\vp_q(\mu_k)=\vt_{q+1}(\mu_k)=0$
and then $J(\mu_k)=0$ Thus the function $\cF$ has at least double
zero at $\mu_k.$

The rest of the proof is similar to the proof of Lemma 3.3 in \cite{K1}. \qed

\section{Small perturbations}

Here we redefine the perturbation coefficients and  write $t v_j,$ $t u_j,$ instead of $ v_j,$ $ u_j.$

For $\l\in\gamma_n=(\l_n^-,\l_n^+)$ we have $1-\Delta^2(\l)<0$ and $\cF(\l)=4(1-\Delta^2)(1+A)^2+J^2<4(1-\Delta^2)+J(t)^2.$
As $J(t)={\mathcal O}(t)$ then for $t$ small enough we have $ \cF(\l) <0.$ As on each zone $(\l_{n-1}^+,\l_n^-)$ we have  $\cF >0$ then there are at least two states on $[\l_n^-,\l_n^+].$

From Lemma \ref{l-3.2} it follows that for $t$ small we have with $\l=\l(t)$ $$\xi(\l)=2(-1)^{q-k+1}\sinh v(\l)(1+{\mathcal O}(t))-J(\l),\qq J(\l)={\mathcal O}(t),\,\,v(\l)={\mathcal O}(t),$$ where  $v=\Im\vk$ and $\pm v(\l) >0$ for $\l\in\gamma_k^\pm.$    Let $\l_0=\l_k^-(t)$ or $\l_0=\l_k^+(t),$ $t>0.$ Then $\l_0$ is the bound (antibound) state iff $-(-1)^{q-k}J(\l_0) >0$ ($-(-1)^{q-k}J(\l_0) <0$ respectively) and  i) in Theorem \ref{Th_representation} follows.

{\bf  Proof of ii) in Theorem \ref{Th_representation}.}

 We denote $y_k=\vp_k^{(j)}$ a sin-type polynomial for the $j-$shifted problem
\[\lb{shifted}a_{k+j-1}^0 y_{k-1}+a_{k+j}^0y_{k+1}+b_{k+j}^0y_k=\l
 y_k,\,\,k\in\Z,\,\,y_0^{(j)}=0,\,\,y_1^{(j)}=1.\]

We have (see Toda \cite{To}) $$\vp_q\psi_n^+\psi^-_n=\frac{a_0^0}{a_n^0}\vp_q^{(n)}$$ and put  $\vp(\l,n+k,k)=\vp_n^{(k)}.$ Then  $\vp(\l,k,k)=0,$  $\vp(\l,k+1,k)=1.$

From \cite{Kh1} or \cite{EMT} it follows that we have representation $$f_n^+=\psi_n^++t\sum_{k=n+1}^{p+1} (\breve{J} K
(\l,n,.))_k f_k^+=\psi_n^++t\sum_{k=n+1}^{p+1} (\breve{J} K
(\l,n,.))_k \psi_k^++{\mathcal O}(t^2),$$ where
\begin{equation}\lb{kernel}
 (\breve{J} K
(\l,n,.))_k=\frac{\vp(\l,n,k-1)}{a_{k-1}^0} u_{k-1}+\frac{\vp(\l,n,k)}{a_{k}^0} v_k+\frac{\vp(\l,n,k+1)}{a_{k+1}^0} u_k.
\end{equation} Moreover
\begin{align*}&\tilde{\vp}_0=\vp_0+t\sum_{k=1}^{p+1} (\breve{J} K
(\l,0,.))_k \vp_k+{\mathcal O}(t^2),\,\,
\tilde{\vp}_1=\vp_1+t\sum_{k=2}^{p+1} (\breve{J} K (\l,1,.))_k
\vp_k+{\mathcal O}(t^2),\\
&\tilde{\vt}_0=\vt_0+t\sum_{k=1}^{p+1} (\breve{J} K (\l,0,.))_k
\vt_k+{\mathcal O}(t^2),\,\,
\tilde{\vt}_1=\vt_1+t\sum_{k=2}^{p+1} (\breve{J} K (\l,1,.))_k
\vt_k+{\mathcal O}(t^2).
\end{align*}
Denote $\sum_{k=n+1}^{p+1} (\breve{J} K (\l,n,.))_k\vp_k=P_n$ and
$\sum_{k=n+1}^{p+1} (\breve{J} K (\l,n,.))_k\vt_k=T_n.$ Then we get
\begin{align*}
&J=-t\left[\vp_qT_1+\phi(P_1-T_0)+
\frac{ v_0}{a_0^0}\vp_q \vt_0+\vt_{q+1}P_0 +\frac{ u_0}{a_0^0}\phi\right]+{\mathcal O}(t^2),\\
&A=  t\frac12\left[
P_1+T_0+\frac{ u_0}{a_0^0}\right]+{\mathcal
O}(t^2) .
\end{align*}
  We denote $J_1$ respectively  $A_1$ the coefficient for $t$ in the expansion:   $J=tJ_1+{\mathcal O}(t^2)$ respectively  $A=tA_1+{\mathcal O}(t^2).$ Then
\begin{align*}
J_1&=-\left[\vp_qT_1+\phi(P_1-T_0)+
\frac{ v_0}{a_0^0}\vp_q +
\vt_{q+1}P_0 + \frac{ u_0}{a_0^0}\phi\right]=-\left[\vp_q\sum_{k=2}^{p+1}\breve{J}K(\l,1,k)\vt_k\right.\\
&\left.+\phi\left(
\sum_{k=2}^{p+1}\breve{J}K(\l,1,k)\vp_k-\sum_{k=2}^{p+1}\breve{J}K(\l,0,k)\vt_k\right)+
\frac{ v_0}{a_0^0}\vp_q+\vt_{q+1}\sum_{k=1}^{p+1}
\breve{J}K(\l,0,k)\vp_k+\frac{ u_0}{a_0^0}\phi\right],\\
A_1&=\frac{1}{2}\left(\sum_{k=2}^{p+1}\breve{J}K(\l,1,k)\vp_k+\sum_{k=1}^{p+1}\breve{J}K(\l,0,k)\vt_k+\frac{ u_0}{a_0^0}\right).
\end{align*}

We restrict ourself to the case $ u_j=0$ as otherwise the formulas would be too cumbersome due to three terms in (\ref{kernel})).  Using that $\vp(\l,0,k)=-\frac{a_k^0}{a_0^0}\vp_k,$ $\vp(\l,1,k)=\frac{a_k^0}{a_0^0}\vt_k$  we get
\begin{align*}
J_1&=-\left[\vp_q\sum_{k=2}^{p}\frac{ v_k}{a_k^0} \vp
(\l,1,k)\vt_k +
\frac{ v_0}{a_0^0}\vp_q+\phi\left(\sum_{k=2}^{p}\frac{ v_k}{a_k^0} \vp
(\l,1,k)\vp_k-\sum_{k=2}^{p}\frac{ v_k}{a_k^0} \vp
(\l,0,k)\vt_k\right)+\right.\\
&+\left.\vt_{q+1}\sum_{k=1}^{p}\frac{ v_k}{a_k^0} \vp
(\l,0,k)\vp_k\right]=\\
&=-\sum_{k=2}^{p}\frac{ v_k}{a_k^0}\left[\vp_q \vp
(\l,1,k)\vt_k +\phi\left( \vp
(\l,1,k)\vp_k- \vp
(\l,0,k)\vt_k\right)+\vt_{q+1} \vp
(\l,0,k)\vp_k\right]-\\&-\vt_{q+1}\frac{ v_1}{a_1^0} \vp
(\l,0,1)\vp_1-
\frac{ v_0}{a_0^0}\vp_q
=-\frac{1}{a_0^0}\sum_{k=2}^p v_k\left[\vp_q\vt_k^2+2\phi\vt_k\vp_k-\vt_{q+1}\vp_k^2
\right]+\frac{ v_1}{a_0^0}\vt_{q+1} -
\frac{ v_0}{a_0^0}\vp_q=\\
&=
-\frac{1}{a_0^0}\left(\sum_{k=2}^p v_kF_k^0\right)+\frac{ v_1}{a_0^0}\vt_{q+1} -
\frac{ v_0}{a_0^0}\vp_q,
\end{align*}
where $F_k^0=\vp_q\vt_k^2+2\phi\vt_k\vp_k-\vt_{q+1}\vp_k^2=\vp_q\psi_k^+\psi_k^-.$
Similar
\begin{align*}
A_1&
=\frac{1}{2}\sum_{k=2}^{p}\frac{ v_k}{a_k^0}\left( \vp
(\l,1,k)\vp_k+\vp
(\l,0,k)\vt_k\right)=\frac{1}{2}\sum_{k=2}^{p}\frac{ v_k}{a_0^0}\left( \vt_k\vp_k-\vp_k
\vt_k\right)=0.
\end{align*}
Let $\l_0=\l_k^-(t)$  and $\l^{(0)}$ be the left endpoint  $\l_k^-$ or let $\l_0=\l_k^+(t)$  and $\l^{(0)}$ be  the right endpoint $\l_k^+.$ Let $\l\equiv\l(t)=\l^{(0)}+t\l^{(1)}+t^2\l^{(2)}+{\mathcal O}(t^3)$ be solution of $\cF(\l)=0.$ Let $J(\l_0)=tJ_1(\l_0)+{\mathcal O}(t^2).$
Then $\l^{(1)}=0$ and
$$\l_0(t)=\l^{(0)}+t^2\frac{J_1^2(\l^{(0)})}{4(\Delta^2)'(\l^{(0)})}+{\mathcal O}(t^3),$$
which can be obtained as in the example in Section \ref{s-special-case}.
Using the properties of the polynomial $\cF_0=4(1-\Delta^2)$  we get $-\cF_0'(\l_k^-)=4(\Delta^2)'(\l_k^-) >0$ and
$-\cF_0'(\l_k^+)=4(\Delta^2)'(\l_k^+)  <0,$ $k=1,\ldots,q-1.$

Using $F_k^0=\vp_q\vt_k^2+2\phi\vt_k\vp_k-\vt_{q+1}\vp_k^2=\vp_q\left(\vt_k+\frac{\phi}{\vp_q}\vp_k\right)^2+\frac{1-\Delta^2}{\vp_q}\vp_k^2,$
we get
\begin{align*}
&J_1(\l_n^\pm)=-\frac{\vp_q(\l_n^\pm)}{a_0^0}\sum_{k=0}^p v_k|\psi_k^+(\l_n^\pm)|^2=-\frac{\vp_q(\l_n^\pm)}{a_0^0}\sum_{k=0}^p
 v_k\left(\vt_k(\l_n^\pm)+\frac{\phi(\l_n^\pm)}{\vp_q(\l_n^\pm)}\vp_k(\l_n^\pm)\right)^2,\\
&J_1(\l)=-\frac{\vp_q}{a_0^0}\sum_{n=0}^p v_n\psi_n^+\psi_n^-=
-\frac{1}{a_0^0}\sum_{n=0}^p v_n F_n^0 =
-\sum_{n=0}^p\frac{ v_n}{a_n^0}\vp_q^{(n)},
\end{align*}
where $(\phi(\l_n^\pm))^2=-\vartheta_{q+1}(\l_n^\pm)\varphi_q(\l_n^\pm),$ which achieves the proof of ii) in Theorem \ref{Th_representation}.  \qed

The Remark after Theorem  \ref{Th_representation} follows from the following inequalities.

If $\l_n^\pm\neq\mu_n,$ $n=1,\ldots, q-1,$  then we have $\overline{\psi^+(\l_n^\pm)}=\psi^-(\l_n^\pm)$ and
 \begin{align*}
&(-1)^{q-n}\vp_q(\l_n^-)=(-1)^{q-n}\frac{a_0^0}{A}\prod_{j=1}^{q-1}
\left(\l_n^--\mu_j\right)>0,\\ &(-1)^{q-1-n}\vp_q(\l_n^+)=(-1)^{q-1-n}\frac{a_0^0}{A}\prod_{j=1}^{q-1}
\left(\l_n^+-\mu_j\right)>0,\\
&(-1)^{q-1}\vp_q(\l_0^-) >0,\qq (-1)^{q-1}\vp_q(\l_q^+) >0.
\end{align*}
  If $ v_k >0,$ $k=1,\ldots, p.$ Then for $n=1,\ldots, q-1,$ we have $(-1)^{q-n+1}J_1(\l_n^-) >0$ and  $\l_n^-(t)$ is the bound state and $\l_n^+(t)$ is the antibound state.

If $\l_n^-=\mu_n,$  for some $n=1,\ldots, q-1,$ then
$$J_1(\l_n^-)=\frac{\vt_{q+1}(\mu_n)}{a_0^0}\sum_{k=1}^{p} v_k(\vp_k(\mu_n))^2=
-\frac{1}{A}\prod_{j=1}^{q-1}(\l_n^--\nu_j)\sum_{k=1}^{p} v_k\vp_k(\mu_n))^2$$ and
$$(-1)^{q-n+1}\vt_{q+1}(\mu_n)=(-1)^{q-n+1}\frac{1}{A}\prod_{j=1}^{q-1}
\left(\l_n^--\nu_j\right)>0.$$ Thus in this case the bound and the antibound states are just swaped with respect to the case $\l_n^\pm\neq\mu_n.$ The proof of the Remark is finished.

\section{Asymptotics}
\begin{theorem}[Asymptotics] \lb{Th-asymptotics-number-of-states}
The function $\xi=2i\Omega\alpha,$ which is analytic on $\Lambda$
and has branch points at $\{E_j\}_{j=0}^{2q-1},$ has the following
asymptotics as $\l\in\gamma^+_q$ and $\l\rightarrow\infty:$
\begin{align*}
\xi(\l)=&\frac{\l^q\prod_{j=0}^pa_j^0}{A
\prod_{j=0}^p a_j}\left[1 +{\mathcal
O}\left(\frac{1}{\l}\right)\right]\,\,\mbox{if}\,\,\l\in\Lambda_+\qq A=\prod_{j=0}^{q-1}a_j^0\\
&=\frac{\l^{2p+q-1}}{A\prod_{j=1}^{p}a_j^0\prod_{j=1}^p a_j}\left[\frac{ v_0}{a_0^0 a_0}
((a_p^0)^2- a_p^2)+{\mathcal O}\left(\frac{1}{\l}\right)\right]
\,\,\mbox{for}\,\,\l\in\Lambda_-\,\,\mbox{and for}\,\,a_p^0\neq a_p\\
&=\frac{\l^{2p+q-2}}{A
\prod_{j=1}^{p}a_j^0\prod_{j=1}^p a_j}\left[\frac{ v_0}{a_0^0 a_0}
(-(a_p^0)^2 v_p)+{\mathcal O}\left(\frac{1}{\l}\right)\right]
\,\,\mbox{for}\,\,\l\in\Lambda_-\,\,\mbox{and
for}\,\,a_p^0= a_p.
\end{align*}

The polynomial $\cF(\l)=\xi(\l)\xi^*(\l)$ has the following asymptotics as $\l\rightarrow\infty:$
\begin{align*}
\cF(\l)&=\frac{\l^{2p+2q-1}}{A^2\prod_{j=0}^{p}a_j^0\prod_{j=0}^p a_j}\left( v_0
((a_p^0)^2- a_p^2)+{\mathcal
O}\left(\frac{1}{\l}\right)\right)\qq\mbox{if}\,\,a_p^0\neq a_p\\
&=
\frac{\l^{2p+2q-2}}{A^2\prod_{j=0}^{p}a_j^0\prod_{j=0}^p a_j}\left( v_0
(-(a_p^0)^2 v_p)+{\mathcal O}\left(\frac{1}{\l}\right)\right)
\qq\mbox{if}\,\,a_p^0= a_p.
\end{align*}
\end{theorem}
The theorem implies that the total number of states of $H$ (counting with multiplicities) is either $2p+2q-1$ states if
$a_p^0\neq a_p$ or    $2p+2q-2$ states if $a_p^0= a_p$ and
$b_p^0\neq b_p.$\\
{\bf Proof.}
The proof follows from the asymptotics of $f_{p-n}^+(\l)$ for $\l\in\Lambda_\pm$
and $\l\rightarrow\infty.$ The asymptotics for $\l\in\Lambda_+$ are well known and can be found for example in \cite{T} and \cite{EMT}. Note that the Jost function $f^+_0(\l)$ on $\Lambda_-$ can be formally obtained  as $(f^+_0(\l))^*$ for $\l\in\Gamma\simeq\Lambda_+$ by iteration of (\ref{pert}) starting with $n=p+1$ and $(f_{p+1}^+(\l))^*=\psi_{p+1}^-(\l),$ $\l\in\Gamma,$ where $\phantom{9}^*$ denotes the complex conjugate.

Let $f_n,$  $\psi_n$  denote either $f_n^+$ respectively $\psi^+_n=\vt_n+m_+\vp_n$ or $(f^+_n)^*$ respectively $\psi^-_n=\vt_n+m_-\vp_n,$ for $\l\in\Gamma.$
 We start with
 $$f_{p+1}=\psi_{p+1},\qq
f_p=\frac{a^0_p}{ a_p}\psi_p.$$
  Put
$\Phi(j)=\frac{\psi_{j+1}}{\psi_j}.$ Thus $\Phi(0)=m_\pm.$ Then
$$\psi_p={\prod_{j=0}^{p-1}}{}^*\Phi(j)=\left\{
                                          \begin{array}{ll}
                                            {\prod_{j=0}^{p-1}}{}\Phi(j) & \mbox{for}\,\,p>0 \\
                                            1 & \mbox{for}\,\, p=0 \\
                                            {\prod_{j=0}^{p-1}}{}(\Phi(j))^{-1} &
\mbox{for}\,\, p<0.
                                          \end{array}
                                        \right.
$$
We have (see Teschl \cite{T}, $a^0(n)\equiv a_n^0,$ $b^0(n)\equiv b_n^0$)
$$\Phi^\pm(\l,n)=\left(\frac{a^0(n)}{\l}\right)^{\pm
1}\left(1\pm\frac{b^0(n+\!\!{\scriptsize
                                               \begin{array}{c}
                                                 1 \\
                                                 0 \\
                                               \end{array}}\!\!)}{\l}
+{\mathcal O}\left(\frac{1}{\l^2}\right)
 \right),\qqq\l\rightarrow\infty.
                                             $$
Put $\Psi(n)=\Phi^{-1}(n).$
$$\Psi^\pm(\l,n)=\left(\frac{a^0(n)}{\l}\right)^{\mp
1}\left(1\mp\frac{b^0(n+\!\!{\scriptsize
                                               \begin{array}{c}
                                                 1 \\
                                                 0 \\
                                               \end{array}}\!\!)}{\l}
+{\mathcal O}\left(\frac{1}{\l^2}\right)
 \right),\qqq\l\rightarrow\infty.
                                             $$
\begin{align*}f_{p-1}&=\frac{(\l- b_p)a_p^0\psi_p- a_p^2\psi_{p+1}}{ a_p a_{p-1}}=
\frac{\psi_{p+1}}{ a_p a_{p-1}}\left(
(\l- b_p)a_p^0\Psi(p)- a_p^2\right)=\circledast;\\
f_{p-2}&=\frac{(\l- b_{p-1}) a_{p-1}\circledast-
 a_{p-1}^2\frac{a_p^0}{ a_p}\psi_{p}}{ a_{p-1} a_{p-2}}=\\
&=\frac{\psi_{p+1}}{ a_p a_{p-1} a_{p-2}}\left((\l- b_{p-1})\left[(\l- b_p)a_p^0\Psi(p)- a_p^2\right]
- a^2_{p-1} a_p^0\Psi(p)\right)=\circledast;\\
f_{p-3}&=\frac{(\l- b_{p-2}) a_{p-2}\circledast-
 a_{p-2}^2\frac{\psi_{p+1}}{ a_p a_{p-1}}\left(
(\l- b_p)a_p^0\Psi(p)- a_p^2\right)}{ a_{p-2} a_{p-3}}=\\
&=\frac{\psi_{p+1}}{ a_p\ldots a_{p-3}}\left(
(\l- b_{p-2})\left[(\l- b_{p-1})\left[(\l- b_p)a_p^0\Psi(p)- a_p^2\right]
- a^2_{p-1} a_p^0\Psi(p)\right]-\right.\\
&-\left. a_{p-2}^2\left(
(\l- b_p)a_p^0\Psi(p)- a_p^2\right)\right).
\end{align*}
Let $\displaystyle
\Psi(p)\equiv\Psi^-(\l,p)=\frac{a_p}{\l}\left(1+\frac{b_p}{\l}
+{\mathcal O}\left(\frac{1}{\l^2}\right)
 \right),$ $\l\rightarrow\infty.$
$$\psi_{p+1}\equiv\psi_{p+1}^-(\l)=\frac{\l^{p+1}}{\prod_{j=0}^{p}a_j^0}\left(1-\frac{1}{\l}\sum_{j=0}^{p}b_j^0
+{\mathcal O}\left(\frac{1}{\l^2}\right)
 \right),\qqq\l\rightarrow\infty.
                                             $$

Then
$$(\l- b_p)a_p^0\Psi(p)- a_p^2=((a_p^0)^2- a_p^2)+\frac{(a_p^0)^2}{\l}((b_p^0)- b_p)+{\mathcal O}\left(\frac{1}{\l^2}\right).$$

\begin{align*}f_{p-n} &=\frac{\l^{p+1}}{\prod_{j=p-n}^p a_j\prod_{j=0}^{p}a_j^0}\\
&\cdot\left( \l^{n-1} ((a_p^0)^2- a_p^2)+\l^{n-2}\left[
-((a_p^0)^2- a_p^2)(\sum_{j=0}^{p}b_j^0+\sum_{j=p-n+1}^{p-1} b_j)+
(a_p^0)^2(b_p^0- b_p)\right]\right)+\\
&+{\mathcal O}\left(\l^{p+n-2}\right);\\
&f_0(\l)=\frac{\l^{2p}}{\prod_{j=0}^p a_j\prod_{j=0}^{p}a_j^0}\\
&\cdot\left( ((a_p^0)^2- a_p^2)+\l^{-1}\left[
-((a_p^0)^2- a_p^2)(\sum_{j=0}^{p}b_j^0+\sum_{j=1}^{p-1} b_j)+
(a_p^0)^2(b_p^0- b_p)\right]\right)+ {\mathcal
O}\left(\l^{2p-2}\right).\end{align*} If $ a_p=a_p^0,$ then
$f_0(\l)=A^{-2}\l^{2p-1}(b_p- b_p^0)+{\mathcal
O}\left(\l^{2p-2}\right).$

Now we can get the asymptotics of  function $\xi.$ Recall that
$$\xi=2i\sin q\vk (\l)
\alpha
(\l)=\frac{\vp_q}{a_0^0}w=\frac{\vp_q}{a_0^0}\left\{f^-,f^+\right\}=\varphi_q(\l)(\frac{ a_0}{a_0^0}f_1^++(\frac{ v_0}{a_0^0}-m_-)f_0^+).$$
Using $m_\pm=\Phi (0)$ and
\begin{align*}
\vp_q&=\frac{\l^{q-1}}{\prod_{j=1}^{q-1} a_j^0}+{\mathcal
O}(\l^{q-2}),\qq m_\pm =\left(\frac{a_0^0}{\l}\right)^{\pm
1}\left(1\pm\frac{b^0(\!\!{\scriptsize
                                               \begin{array}{c}
                                                 1 \\
                                                 0 \\
                                               \end{array}}\!\!)}{\l}
+{\mathcal O}\left(\frac{1}{\l^2}\right)
 \right),\\
f_0^+&=\alpha_0^+\left[1 +{\mathcal
O}\left(\frac{1}{\l}\right)\right],\qq
f_1^+=\alpha_1^+\frac{a_0^0}{\l}\left[1
+{\mathcal O}\left(\frac{1}{\l}\right)\right],
\end{align*}
\begin{align*}
(f_0^+)^*&=\frac{\l^{2p}}{\prod_{j=0}^p a_j\prod_{j=0}^{p}a_j^0}\\
&\cdot\left(  ((a_p^0)^2- a_p^2)+\l^{-1}\left[
-((a_p^0)^2- a_p^2)(\sum_{j=0}^{p}b_j^0+\sum_{j=1}^{p-1} b_j)+
(a_p^0)^2(b_p^0- b_p)\right]\right)+{\mathcal
O}\left(\l^{2p-2}\right),
\end{align*}
\begin{align*}
(f_1^+)^*&=\frac{\l^{2p-1}}{\prod_{j=1}^p a_j\prod_{j=0}^{p}a_j^0}\\
&\cdot\left(  ((a_p^0)^2- a_p^2)+\l^{-1}\left[
-((a_p^0)^2- a_p^2)(\sum_{j=0}^{p}b_j^0+\sum_{j=2}^{p-1} b_j)+
(a_p^0)^2(b_p^0- b_p)\right]\right)+{\mathcal
O}\left(\l^{2p-3}\right),
\end{align*} where
$$ \alpha_0^+=\prod_{k=0}^p\frac{a_k^0}{ a_k},\,\,
\alpha_1^+=\prod_{k=1}^p\frac{a_k^0}{ a_k},\,\,
\alpha_0^-=\prod_{k=-1}^0\frac{a_k^0}{ a_k},\,\,
\alpha_1^-=\frac{a_0^0}{ a_0},
$$
we get, for $\l\in\Lambda_+,$
\begin{align*}
\xi&=\frac{\l^{q-1}}{\prod_{j=1}^{q-1} a_j^0}\left[
\frac{\alpha_1^+ a_0}{\l}\left\{1 +{\mathcal
O}\left(\frac{1}{\l}\right)\right\} + \left(\frac{\l}{a_0^0}+{\mathcal
O}(1)\right)\alpha_0^+\left\{1 +{\mathcal
O}\left(\frac{1}{\l}\right)\right\}\right],\\
&=\frac{\l^q\alpha_0^+}{a_0^0\prod_{j=1}^{q-1} a_j^0}\left[1 +{\mathcal
O}\left(\frac{1}{\l}\right)\right],
\end{align*}
and for $\l\in\Lambda_-,$
\begin{align*}
\xi&=\frac{\l^{q-1}}{\prod_{j=1}^{q-1}
a_j^0\prod_{j=0}^{p}a_j^0\prod_{j=1}^p a_j}\left[\frac{ a_0}{a_0^0}
\l^{2p-1} \left\{ ((a_p^0)^2- a_p^2) +{\mathcal
O}\left(\frac{1}{\l}\right)\right\}+\right.\\
&\left.+\left(\frac{ v_0}{a_0^0}+{\mathcal
O}\left(\frac{1}{\l}\right)\right)\frac{\l^{2p}}{ a_0}\left\{
((a_p^0)^2- a_p^2)+{\mathcal
O}\left(\frac{1}{\l}\right)\right\}\right]=\\
&=\frac{\l^{2p+q-1}}{\prod_{j=1}^{q-1}
a_j^0\prod_{j=0}^{p}a_j^0\prod_{j=1}^p a_j}\left[\frac{ v_0}{a_0^0 a_0}
((a_p^0)^2- a_p^2)+{\mathcal O}\left(\frac{1}{\l}\right)\right]
\qq\mbox{or}\\ &=\frac{\l^{2p+q-2}}{\prod_{j=1}^{q-1}
a_j^0\prod_{j=0}^{p}a_j^0\prod_{j=1}^p a_j}\left[\frac{ v_0}{a_0^0 a_0}
(-(a_p^0)^2 v_p)+{\mathcal
O}\left(\frac{1}{\l}\right)\right]\qq\mbox{if}\,\,a_p= a_p.
 \end{align*}
 Hence the polynomial $\cF(\l)=\xi(\l)\xi^*(\l)$
has asymptotics
\begin{align*}
\cF(\l)&=\frac{\l^{2p+2q-1}}{A^2\prod_{j=0}^{p}a_j^0\prod_{j=0}^p a_j}\left( v_0
((a_p^0)^2- a_p^2)+{\mathcal
O}\left(\frac{1}{\l}\right)\right)\qq\mbox{or}\\
&=
\frac{\l^{2p+2q-2}}{A^2\prod_{j=0}^{p}a_j^0\prod_{j=0}^p a_j}\left( v_0
(-(a_p^0)^2 v_p)+{\mathcal O}\left(\frac{1}{\l}\right)\right)
\qq\mbox{if}\,\,a_p^0= a_p.
\end{align*}
 \qed

\section{Example: $p=1,$ $ a_j=a_j^0.$}\lb{s-special-case}
Let $p=1,$ $ a_j=a_j^0.$ Then $f_1^+=\psi_1^+=m_+.$ Using that
$$\psi_2^+=\frac{(\l-b_1^0)\psi_1^+-a_0^0\psi_0^+}{a_1^0}=\frac{(\l-b_1^0)m_+-a_0^0}{a_1^0}$$ we get
$$f_0^+=\frac{(\l- b_1)\psi_1^+-a_1^0\psi_2^+}{a_0^0}=\frac{(\l - b_1)m_+-(\l-b_1^0)m_+
+a_0^0}{a_0^0}=1-\frac{ v_1}{a_0^0} m_+.$$ Thus we have
$\tilde{\vt}_0=1,$ $\tilde{\vp}_0=-\frac{ v_1}{a_0^0},$
$\tilde{\vt}_1=0,$ $\tilde{\vp}_1=1.$ Then
$$\xi(\l)=\frac{ v_0}{a_0^0} \vp_q(\l)
-\frac{ v_1}{a_0^0}\vt_{q+1}+(1-\frac{ v_0 v_1}{(a_0^0)^2})m_+\vp_q(\l)
-m_- \vp_q(\l).$$ Now $$J=-\frac{ v_0}{a_0^0}\vp_q
+\frac{ v_1}{a_0^0}\vt_{q+1}+\frac{ v_0 v_1}{(a_0^0)^2}\phi,\qq
A=-\frac{ v_0 v_1}{2(a_0^0)^2}.$$ If we write
$t v_j$ instead of $ v,$ then in the first order in
$t$ we get
$$\cF =4(1-\Delta^2)(1+A)^2+J^2=4(1-\Delta^2)+t^2\left(-\frac{ v_0}{a_0^0}\vp_q
+\frac{ v_1}{a_0^0}\vt_{q+1}\right)^2+{\mathcal O}(t^3).$$
Denote $(-\frac{ v_0}{a_0^0}\vp_q
+\frac{ v_1}{a_0^0}\vt_{q+1})^2=\cF_2(\l).$
 Let
$\l^{(0)}=\l_k^+$ or $\l^{(0)}=\l_k^-$ for some $k=1,\ldots, q-1.$ We look for
solutions of $\cF(\l)=0$  in the form $\l(t)=\l^{(0)}+\l^{(1)}t +
\l^{(2)} t^2+\ldots.$ As $4(1-\Delta^2(\l^{(0)}))=0$ we get the Taylor expansion at
$\l=\l^{(0)}:$ \begin{align*}\cF(\l)=&-4(\Delta^2)'(\l^{(0)})(\l^{(1)}t +
\l^{(2)}t^2+\ldots)  -2(\Delta^2)''(\l^{(0)})(\l^{(1)}t +
\l^{(2)}t^2+\ldots)^2+\ldots\\
&+t^2\cF_2(\l^{(0)})+t^2\cF_2(\l^{(0)})'(\l^{(1)}t + \l^{(2)}t^2+\ldots)+{\mathcal
O}(t^3).\end{align*}

As $(\Delta^2)'(\l^{(0)})\neq 0,$ $(\Delta^2)''(\l^{(0)})\neq 0,$ we get
$\l^{(1)}=0,$
$$\l^{(2)}=\frac{\cF_2(\l^{(0)})}{4(\Delta^2)'(\l^{(0)})}=\frac{1}{4(\Delta^2)'(\l^{(0)})}
\left(-\frac{ v_0}{a_0^0}\vp_q(\l^{(0)})
+\frac{ v_1}{a_0^0}\vt_{q+1}(\l^{(0)})\right)^2$$ and
$\l(t)=\l^{(0)}+\l^{(2)}t^2+{\mathcal O}(t^3).$

\end{document}